\documentclass[a4paper, 11pt]{article}
\usepackage[english]{babel}
\DeclareMathAlphabet\mathbfcal{OMS}{cmsy}{b}{n}
\usepackage{latexsym}
\usepackage{amsfonts, amsthm}
\usepackage{amssymb}
\usepackage{amscd}
\usepackage{amsmath}
\usepackage{graphics}
\usepackage{eepic}
\usepackage{epsfig}
\usepackage[noxcolor]{pstricks}
\newtheorem{theorem}{Theorem}
\newtheorem{lemma}[theorem]{Lemma}
\newtheorem{proposition}[theorem]{Proposition}
\newtheorem{corollary}[theorem]{Corollary}
\newtheorem{definition}
{Definition}
\newtheorem{remark}{Remark\/}
\newtheorem{example}{Example\/}
{Conjecture\/}

\newcommand{\C}{{\mathbb C}}
\newcommand{\Z}{{\mathbb Z}}

\newcommand{\I}{\mathcal{I}}

\begin{document}

\title{\textbf{Poincar\'e series for filtrations defined by discrete valuations with arbitrary center}}
\author{Antonio Campillo and Ann Lemahieu
\date{}}
\maketitle \noindent {\footnotesize \emph{\noindent \textbf{Abstract.---} To study singularities on complex varieties we study Poincar\'e series of filtrations that are defined by discrete valuations on the local ring at the singularity. In all previous papers on this topic one poses restrictions on the centers of these valuations and often one uses several definitions for Poincar\'e series. In this article we show that these definitions can differ when the centers of the valuations are not zero-dimensional, i.e.\ do not have the maximal ideal as center.
We give a unifying definition for Poincar\'e series which also allows filtrations defined by valuations that are all nonzero-dimensional. We then show that this definition satisfies a nice relation between Poincar\'e series for embedded filtrations and Poincar\'e series for the ambient space and we give some application for singularities which are nondegenerate with respect to their Newton polyhedron.}}
\\ \\
${}$
\begin{center}
\textsc{0. Introduction}
\end{center}
${}$
\\
In \cite{cdk} one introduced a Poincar\'e series induced by a filtration on the ring of germs of a complex variety.
This Poincar\'e series has been studied for several kinds of singularities, see for example \cite{triomonodromy}, \cite{trio1}, \cite{trio3}, \cite{ebeling}, \cite{cutkosky}, \cite{lemahieu}, \cite{pedro} and \cite{nemethi}. In some cases this Poincar\'e series determines the topology of the singularity and is related to its zeta function of monodromy.
\\
\indent In some of these papers one uses several definitions for Poincar\'e series. In Section 1 we will see that these definitions become not necessarily equivalent when not all considered valuations are centered at the maximal ideal. In \cite{pedro}, one studies the Poincar\'e series for quasi-ordinary and for toric singularities and, in particular, a new definition for Poincar\'e series was introduced for this context which made it possible to treat also sets of monomial valuations where at least one valuation was centered at the maximal ideal.
In Section 2 we give a unifying definition for the Poincar\'e series, i.e.\ a definition that coincides with the former definitions when all valuations are centered at the maximal ideal, that also coincides with the one in \cite{pedro} when the valuations are monomial and at least one of them is centered at the maximal ideal, and that even makes sense in some cases where none of the valuations is centered at the maximal ideal. This Poincar\'e series is defined in homological terms.
\\ \indent
In \cite{annembedded} the second author introduced a Poincar\'e series for embedded varieties in an ambient space by taking multi-index filtrations coming from valuations on the ambient space. When at least one such valuation is centered at the maximal ideal and when the subspace corresponds to a principal ideal, then a nice formula relating the Poincar\'e series of the embedded space and the ambient space showed up. It was shown that this relating formula gives rise to interesting topological and geometrical information in the case of plane curve singularities and singularities that are nondegenerate with respect to their Newton polyhedron.
We can now also define Poincar\'e series for embedded filtrations in homological terms, see Section 3. We show that the relating formula from \cite{annembedded} still holds when none of the valuations is centered at the maximal ideal, under the condition that both Poincar\'e series are well defined. We then extend the topological and geometrical results from \cite{annembedded} to the context of sets of valuations where none of them is centered at the maximal ideal.

${}$
\begin{center}
\textsc{1. Poincar\'e series - different definitions}
\end{center}
${}$ \\
In this section we would like to call attention to the fact that
Poincar\'e series have been defined in several ways, at first sight maybe
equivalent up to notation, but in fact different in the sense of less or more general.
${}$ \\ \\
\emph{\textbf{Description 1.--}} Let $(X,o)$ be a germ of a complex algebraic variety and let
$\mathcal{O}_{X,o}$ be the local ring of germs of functions on
$(X,o)$. Let $\underline{\nu}=\{\nu_1,\cdots,\nu_r\}$ be a set of
order functions from $\mathcal{O}_{X,o}$ to $\mathbb{Z} \cup
\{\infty\}$, i.e.\ functions $\nu_j$ that satisfy $\nu_j(f+g) \geq
\textnormal{ min } \{\nu_j(f),\nu_j(g)\}$ and $\nu_j(fg) \geq
\nu_j(f)$, for all $f, g \in \mathcal{O}_{X,o}$, $1 \leq j \leq r$.
An order function $\nu_j$ on $\mathcal{O}_{X,o}$ is called a valuation if moreover it
satisfies $\nu_j(fg) = \nu_j(f)+\nu_j(g)$, for all $f, g \in
\mathcal{O}_{X,o}$.
 The set
$\underline{\nu}$ defines a multi-index filtration on
$\mathcal{O}_{X,o}$ by the ideals
\[M(\underline{v}):=\{g \in \mathcal{O}_{X,o} \mbox{ $|$ } \nu_j(g)
\geq v_j, 1 \leq j \leq r\}, \qquad \underline{v} \in
\mathbb{Z}^r.\] If the dimensions of the complex vector spaces
$M(\underline{v})/M(\underline{v}+\underline{1})$ are finite for all
$\underline{v} \in \mathbb{Z}^r$, then originally (see \cite{cdk}
and \cite{trio1}) the \emph{Poincar\'e series} associated to this
multi-index filtration was defined as
\begin{eqnarray} \label{deforiginal}
P^{\underline{\nu}}_{X}(t_1,\ldots,t_r):=\frac{\prod_{j=1}^{r}(t_j-1)}{(t_1\cdot \ldots \cdot
t_r - 1)}\sum_{\underline{v} \in \mathbb{Z}^r}
\mbox{dim}(M(\underline{v})/M(\underline{v}+\underline{1}))
\underline{t}^{\underline{v}}.
\end{eqnarray}
In particular, one may have the valuative case, i.e. when
$\underline{\nu}$ is a set of
discrete valuations of the function field $\mathbb{C}(X)$ whose
valuation rings contain $\mathcal{O}_{X,o}$. Then $\{f \in \mathcal{O}_{X,o} | \nu_j(f) > 0\}$ is a prime ideal which is called
the center of the valuation $\nu_j$, $1 \leq j \leq r$.
In the valuative case, this description is well defined if the
center of each valuation is the maximal ideal $\mathsf{m}$ of
$\mathcal{O}_{X,o}$ (see \cite{lemahieu}). \\
\indent One interesting valuative case is the toric one, i.e. when the data are an affine toric variety and the $\nu_j$ are monomial valuations. Let $N = \mathbb{Z}^n, n>1$, and let $M$ be the dual space to $N$, then there is a natural bilinear map $M \times N \rightarrow \mathbb{Z} : (m , n) \mapsto \langle m , n \rangle$. The dual cone  $\check{\sigma}$ to a cone  $\sigma \subset N \otimes_{\mathbb{Z}}\mathbb{R}$ is defined as
the set $\{m \in \mathbb{R}^n \mid \langle m,x \rangle \geq 0,
\forall x \in \sigma\}$. Now one considers a semigroup $S$ generating $\check{\sigma}$ as a cone, and the affine toric variety $X=$ Spec $\C[S]$. For each $\nu \in \sigma$, the duality allows to define a discrete valuation on $\C[S]$, also denoted by $\nu$, given by $\nu(x) = \langle \nu,x \rangle$ on the monomials $x$ of $\C[S]$. Such valuations, which are centered on a prime ideal of $\C[S]$, are called monomial, since their value on a function $f \in \C[S]$ is the minimum of the values $v(x)$ for $x$ in the support of $f$.
\\ \\
\emph{\textbf{Description 2.--}}
For $I$ a subset in $\{1,\ldots,r\}$, let $\underline{e}_I$ be the
$r$-tuple with $j$-th component equal to $1$ if $j \in I$ and equal
to $0$ otherwise. When $I$ is a singleton, say $I=\{j\}$, we also denote $\underline{e}_I$ by $\underline{e}_j$.
Notice that the coefficient of $\underline{t}^{\underline{v}}$ in the Poincar\'e series (1) can also be written as
\begin{eqnarray*}
& & \sum_{I \subset
\{1,\ldots,r\}}(-1)^{\#I}\mbox{dim}\frac{M(\underline{v}+\underline{e}_I)}{M(\underline{v}+\underline{1})}\\
& = & \sum_{K \subset
\{1,\ldots,r-1\}}(-1)^{\#K}\mbox{dim}\frac{M(\underline{v}+\underline{e}_K)}{M(\underline{v}+\underline{e}_K+\underline{e}_r)}.
\end{eqnarray*}
However, this last expression is well defined when the vector spaces $\frac{M(\underline{v}+\underline{e}_K)}{M(\underline{v}+\underline{e}_K+\underline{e}_r)}$ are of finite dimension,
what in the valuative case means that only $\nu_r$ needs to be centered at $\mathsf{m}$.
\\ \\
\emph{\textbf{Description 3.--}}
When only considering two valuations, the coefficient of $\underline{t}^{\underline{v}}$ in the Poincar\'e series (1) can also be written as
\[\mbox{dim } \frac{M(v_1,v_2)}{M(v_1+1,v_2)+M(v_1,v_2+1)}.\]
This term can be well defined although none of the valuations is centered
at $\mathsf{m}$, as the following example shows. We will see in Section 2 that this description can be given in a much more general context.
\begin{example}
\emph{\label{excenter}
We take $\sigma$ to be the cone $\mathbb{R}_{\geq 0}^2$.
 Let $S \subset \check{\sigma}
\cap M$ be the semigroup generated by the
vectors $(1,0)$ and $(0,1)$ and let $X$ be the affine toric variety Spec $\C[S]$.
We consider the monomial valuations
$\nu_1$ and $\nu_2$ corresponding to the vectors $(1,0)$ and $(0,1)$
in $\sigma$, i.e.\ $\nu_1(x^a y^b)=a$ and
$\nu_2(x^a y^b)=b$ for $(a,b) \in \mathbb{Z}_{\geq 0}^2$. These valuations are obviously not centered at the maximal ideal.
One computes that \mbox{dim} $M(v_1,v_2)/(M(v_1+1,v_2)+M(v_1,v_2+1))=1$,
for all $(v_1,v_2) \in \mathbb{Z}_{\geq 0}^2$ and that \mbox{dim} $M(v_1,v_2)/(M(v_1+1,v_2)+M(v_1,v_2+1))=0$,
for all $(v_1,v_2) \notin \mathbb{Z}_{\geq 0}^2$.
Hence we get that the Poincar\'e series defined in this way is equal to $\frac{1}{(1-t_1)(1-t_2)}.$
\hfill
\vspace*{-0.5cm}$\square$}
\end{example}
${}$
\\ \\
\emph{\textbf{Description 4.--}}
In the paper \cite{trio1} on Poincar\'e series for plane curve
singularities, an alternative description was given for this
Poincar\'e series: for $1 \leq j \leq r$, denote by
$D_j(\underline{v})$ the complex vector space
$M(\underline{v})/M(\underline{v}+\underline{e}_j)$. Let us consider the
map
\begin{eqnarray*}
j_{\underline{v}}: M(\underline{v}) & \longrightarrow &
D_1(\underline{v}) \times \cdots \times D_r(\underline{v}) \\
g & \longmapsto & (a_1(g),\ldots,a_r(g))=:\underline{a}(g),
\end{eqnarray*}
where $a_j(g)$ is the projection of $g$ on $D_j(\underline{v})$. Let
$D(\underline{v})$ be the image of the map $j_{\underline{v}}$, then
$D(\underline{v})\simeq
M(\underline{v})/M(\underline{v}+\underline{1})$. One defines the
fibre $F_{\underline{v}}$ as the space $D(\underline{v}) \cap
(D_1^{*}(\underline{v}) \times \cdots \times
D_r^{*}(\underline{v}))$, where $D_j^{*}(\underline{v})$ denotes
$D_j(\underline{v}) \setminus \{0\}$. Then
$F_{\underline{v}}$ is invariant with respect to multiplication by
nonzero constants; let
$\mathbb{P}F_{\underline{v}}:=F_{\underline{v}}/\mathbb{C}^{*}$ be
the projectivisation of $F_{\underline{v}}$. Let $\chi$ denote the topological Euler characteristic.
\begin{theorem} (\cite[Theorem 3]{trio1})\label{def2}
If the valuations $\nu_1,\ldots,\nu_r$ are centered at the maximal ideal, then
$P^{\underline{\nu}}_{X}(\underline{t}) = \sum_{\underline{v} \in
\mathbb{Z}^r}\chi(\mathbb{P}F_{\underline{v}})\underline{t}^{\underline{v}}$.
\end{theorem}
\noindent This equality holds for all singularities $(X,o)$ (not
only plane curve singularities), but the proof only works if all valuations are centered at the maximal ideal.
When this is not the case, it can happen that $\chi(\mathbb{P}F_{\underline{v}})$ is infinite but it can also happen that $\chi(\mathbb{P}F_{\underline{v}})$ is finite for all $\underline{v} \in
\mathbb{Z}^r$ and thus that the expression $\sum_{\underline{v} \in
\mathbb{Z}^r}\chi(\mathbb{P}F_{\underline{v}})\underline{t}^{\underline{v}}$
is well defined. In the toric case this expression is easy to compute, as shown in the following example.
\\ \\
\noindent \textbf{Example 1 Continued.}\\
   In the toric monomial case one can decompose $F_{\underline{v}}$ as a (here infinite) disjoint union according to the support of the functions in $F_{\underline{v}}$.
   We see that only the monomials contribute to $\chi(\mathbb{P}F_{(v_1,v_2)})$. In particular,
when $\chi(\mathbb{P}F_{\underline{v}})$ is finite, then it equals the number of monomials with value equal to $\underline{v}$ (see also \cite{lemahieu}, Proof of Prop.\ 1).
In this example we thus find $\chi(\mathbb{P}F_{(v_1,v_2)})=1$ for all $(v_1,v_2)
\in \mathbb{Z}^2_{\geq 0}$ and $\chi(\mathbb{P}F_{(v_1,v_2)})=0$ for all $(v_1,v_2)
\notin \mathbb{Z}^2_{\geq 0}$. Hence
$\sum_{\underline{v} \in \mathbb{Z}^r}
\chi(\mathbb{P}F_{\underline{v}})\underline{t}^{\underline{v}}=\frac{1}{(1-t_1)(1-t_2)}.$  \hfill
\vspace*{-0.5cm} $\square$
\\ \\ \\
In the next section we will give a unifying definition for the above described Poincar\'e series.
\\
${}$
\begin{center}
\textsc{2. Poincar\'e series defined in homological terms}
\end{center}
${}$ \\
Fix some $\underline{v}
\in \mathbb{Z}^r$. We will denote
$V_I:=M(\underline{v}+\underline{e}_I)/M(\underline{v}+\underline{1})$
for $I \subseteq \{1,\ldots,r\}$.
Let \[ C_i:=\begin{cases}
\begin{array}{cc}
\{0\} & \mbox{ if } i=-1, i>r, \\
\bigoplus_{\substack{I
\subset \{1,\ldots,r\}, \\ \#I=i}} V_I & \mbox{ if } 0 \leq i \leq r.
\end{array}
\end{cases} \]
For $-1 \leq i \leq r-1$, we define a map $\partial_{i+1}:
C_{i+1} \rightarrow C_i$ by defining it on each component $V_I$ ($\#
I=i+1$). Suppose that $I=\{a_1,a_2,\ldots,a_{i+1}\}\subseteq \{1,\ldots,r\}$ with $a_1 < a_2 < \ldots < a_{i+1}$, then we set
\begin{eqnarray*}
V_I & \longrightarrow & C_i=\bigoplus_{\substack{J
\subset \{1,\ldots,r\}, \\{\#J=i}}} V_{J} \\
x & \longmapsto & (\underline{y})_J,
\end{eqnarray*}
where $y_J=0$ if $J \nsubseteq I$ and $y_J=(-1)^k x$ if $J = I
\setminus \{a_k\}$. For instance, if $r=4$, then

\footnotesize{
\begin{eqnarray*}
\partial_3: V_{\{1,2,3\}}\oplus V_{\{1,2,4\}}\oplus V_{\{1,3,4\}}\oplus
V_{\{2,3,4\}} & \longrightarrow & V_{\{1,2\}}\oplus
V_{\{1,3\}}\oplus V_{\{1,4\}}\oplus V_{\{2,3\}}\oplus
V_{\{2,4\}}\oplus V_{\{3,4\}} \\
(x,y,z,u) & \longmapsto & (-x-y,x-z,y+z,-x-u,-y+u,-z-u).
\end{eqnarray*}}

\normalsize \noindent For $i \leq r$, we define $\partial_i$ to be the zero map.
Notice that $\partial_i \circ
\partial_{i+1} \equiv 0$ ($0 \leq i \leq r-1$), so $\mathcal{C}_\bullet=(C_i,\partial_i)_{i \in \mathbb{Z}_{\geq -1}}$ defines a complex of vector spaces. Let
$H_i(\mathcal{C})=\textnormal{ Ker } \partial_i / \textnormal{ Im }
\partial_{i+1}$ and $h_i=\textnormal{ dim }H_i(\mathcal{C})$. If confusion about the considered
vector $\underline{v}$ in $\Z^r$ is possible, we will write
$h_i^{\underline{v}}$ and $\mathcal{C}^{\underline{v}}_{\bullet}=
(C^{\underline{v}}_i,\partial^{\underline{v}}_i)$.


\begin{theorem} \label{thmhomologia}
If $\nu_1, \ldots, \nu_r$ are discrete valuations of $X$ and if $\nu_r$ is centered at the maximal ideal, then the Poincar\'e series defined in \emph{Description 1} or \emph{Description 2} coincides with
\[\sum_{\underline{v} \in \mathbb{Z}^r}
\left(\sum_{i = 0}^{r-1}(-1)^i h_i^{\underline{v}} \right)
\underline{t}^{\underline{v}}.\]
\end{theorem}
\noindent \emph{Proof.} \quad As $\nu_r$ has its center at $\mathsf{m}$, the Poincar\'e series of \emph{Description 1} can be rewritten as in \emph{Description 2}, and so it is sufficient to study the Poincar\'e series of \emph{Description 2}.
We first define a complex of finite dimensional vector spaces. Fix again $\underline{v}
\in \mathbb{Z}^r$. For $I \subseteq \{1,\ldots,r-1\}$, let
$\tilde{V}_I:=M(\underline{v}+\underline{e}_I)/M(\underline{v}+\underline{e}_I+\underline{e}_r)$
and let
\[ \tilde{C}_i:=\begin{cases}
\begin{array}{cc}
\{0\} & \mbox{ if } i=-1, i>r-1, \\
\bigoplus_{\substack{I
\subset \{1,\ldots,r-1\}, \\ \#I=i}} \tilde{V}_I & \mbox{ if } 0 \leq i \leq r-1.
\end{array}
\end{cases} \]
For $-1 \leq i \leq r-2$, we define a map $\tilde{\partial}_{i+1}:
\tilde{C}_{i+1} \rightarrow \tilde{C}_i$ by defining it on each component $\tilde{V}_I$ ($\#
I=i+1$), analogously as in the construction of the complex $(C_i,\partial_{i})$, and so we get a new complex of vector spaces $\tilde{\mathcal{C}}_\bullet=(\tilde{C}_i,\tilde{\partial}_i)_{i \in \mathbb{Z}_{\geq -1}}$.
For $i \leq r-1$, we define $\partial_i$ to be the zero map.
Let $\tilde{H}_i(\mathcal{C})=\textnormal{ Ker } \tilde{\partial}_i / \textnormal{ Im }
\tilde{\partial}_{i+1}$ and $\tilde{h}_i=\textnormal{ dim }\tilde{H}_i(\mathcal{C})$. If confusion about the considered
vector $\underline{v}$ in $\Z^r$ is possible, we will write
$\tilde{h}_i^{\underline{v}}$.

For $i \in \{0,\ldots,r-1\}$, consider now the map $\phi_i: C_i \rightarrow \tilde{C}_i$ that sends
$x \in V_I$ to $(\overline{0},\ldots,\overline{0}) \in \tilde{C}_i$ if $r \in I$, and to $(\overline{0},\ldots, \overline{0},\overline{x},\overline{0},\ldots,\overline{0}) \in \tilde{C}_i$ if $r \notin I$,
with $\overline{x}$ on the component of index $\tilde{V}_I$.
We set
\begin{eqnarray*}
L_i &:=& \mbox{ Ker }(\phi_i) \\
& = &  \bigoplus_{\substack{I \in \{1,\ldots,r\}\\r \in I, \# I=i}} V_I \oplus \bigoplus_{\substack{I \in \{1,\ldots,r-1\}\\\# I=i}} \frac{M(\underline{v}+\underline{e}_I+\underline{e}_r)}{M(\underline{v}+\underline{1})}\\
& = & \bigoplus_{\substack{K \in \{1,\ldots,r-1\}\\\# K=i}} \frac{M(\underline{v}+\underline{e}_K+\underline{e}_r)}{M(\underline{v}+\underline{1})}
\oplus \bigoplus_{\substack{J \in \{1,\ldots,r-1\}\\\# J=i-1}} \frac{M(\underline{v}+\underline{e}_J+\underline{e}_r)}{M(\underline{v}+\underline{1})},
\end{eqnarray*}
and we denote the induced maps $\partial^L_i: L_i \rightarrow L_{i-1}$ giving rise to the complex
$\mathcal{L}_\bullet=(L_i,\partial^L_i)_{i \in \mathbb{Z}_{\geq -1}}$, with induced homologies
$h^L_i$.
If $\chi(\mathcal{C}_\bullet)$ is finite, we have $\chi(\mathcal{C}_\bullet)= \chi(\tilde{\mathcal{C}}_\bullet)+\chi(\mathcal{L}_\bullet)$.
We will now prove that $\chi(\mathcal{L}_\bullet)=0$.
Asking that the following diagram commutes
\\${}$
\begin{eqnarray*}
  0   \longrightarrow  L_i  \quad \longrightarrow  C_i  \quad \longrightarrow   \tilde{C}_i  \quad \longrightarrow  0\\
     \downarrow   \partial^L_i  \qquad  \downarrow  \partial_i  \qquad \downarrow  \tilde{\partial}_i  \qquad \quad    \\
0  \longrightarrow  L_{i-1} \longrightarrow  C_{i-1}  \longrightarrow \tilde{C}_{i-1}  \longrightarrow  0
\end{eqnarray*}${}$ \\
means that for $f \in \bigoplus_{\substack{K \in \{1,\ldots,r-1\}\\\# K=i}} \frac{M(\underline{v}+\underline{e}_K+\underline{e}_r)}{M(\underline{v}+\underline{1})}$ and
$g \in \bigoplus_{\substack{J \in \{1,\ldots,r-1\}\\\# J=i-1}} \frac{M(\underline{v}+\underline{e}_J+\underline{e}_r)}{M(\underline{v}+\underline{1})}$,
the differential $\partial^L_i$ acts as $\partial^L_i(f,g)=(\partial_i(f)+(-1)^i g, \partial_{i-1}(g))$.
Now we can deduce that $h^L_i=0$ for all $i$. Indeed, take $(f,g) \in \mbox{ Ker }\partial^L_i$. Then in particular $\partial_i(f)=(-1)^{i+1}g$ and $(f,g)=\partial^L_{i+1}(0,(-1)^{i+1} f)$ what implies that $\chi(\mathcal{L}_\bullet)=0$ and so
$\chi(\mathcal{C}_\bullet)= \chi(\tilde{\mathcal{C}}_\bullet)$.
\\As $\tilde{\mathcal{C}}_\bullet$ is a complex of vector spaces of finite dimension, we can compute $\chi(\tilde{\mathcal{C}}_\bullet)$ as $\sum (-1)^i \mbox{dim} \tilde{C}_i$ what also equals
$\sum_{K \subset
\{1,\ldots,r-1\}}(-1)^{\#K}\mbox{dim}
\frac{M(\underline{v}+\underline{e}_K)}{M(\underline{v}+\underline{e}_K+\underline{e}_r)}$. We can now close the proof.
 \hfill $\blacksquare$
\\ \\
We now compare $\sum_{\underline{v} \in \mathbb{Z}^r}
\left(\sum_{i = 0}^{r-1}(-1)^i h_i^{\underline{v}} \right)
\underline{t}^{\underline{v}}$ with \emph{Description 3} and \emph{Description 4} to motivate a new definition for Poincar\'e series.
\\ \\
\emph{Description 3} goes beyond the domain of sets of valuations of which at least one is centered at the maximal ideal.
Example 1 is illustrating this. Indeed, there one has $h_0^{\underline{v}}= \mbox{dim } \frac{M(v_1,v_2)}{M(v_1+1,v_2)+M(v_1,v_2+1)}$ and $h_1^{\underline{v}}=0$ and so again the Poincar\'e series equals $\sum_{\underline{v} \in \mathbb{Z}^2} (h_0^{\underline{v}}-h_1^{\underline{v}})\underline{t}^{\underline{v}}$.
\\ \\
\emph{Description 4} was only defined in a context where all valuations were centered at the maximal ideal (see
\cite{trio1}), but in Section 1 we have seen how to extend and compute it in the toric monomial case.
\\ \\
With the following proposition, we get a homological description in Corollary \ref{coroltoric}
in the case of affine toric varieties and monomial valuations.

\begin{proposition} \label{propbases}
Suppose that there exist bases $\mathcal{B}_i$ for the vector spaces
$V_{\{i\}}$, $1 \leq i \leq r$, such that $\mathcal{B}=\cup
\mathcal{B}_i$ is a set of linearly independent vectors in
$V_{\emptyset}=M(\underline{v})/M(\underline{v}+\underline{1})$.
Then the complex $\mathcal{C}_{\bullet}=(C_i,\partial_i)$ is exact in $C_i$
for $i \geq 1$, i.e.\ $h_i^{\underline{v}}=0$ for $1 \leq i \leq
r-1$.
\end{proposition}
\noindent \emph{Proof.} \quad
We will prove that $\textnormal{dim (Im}
\partial_i)= \sum_{m \geq 0} (-1)^m \textnormal{ dim }C_{i+m}$.
The statement can then be deduced from this equation. For $I, J \subset \{1,\ldots,r\}$, first notice
that $V_I \cap V_J=V_{I \cup J}$ and that $\cap_{i \in I}
\mathcal{B}_i$ is a basis of $V_I$, as $\mathcal{B}=\cup
\mathcal{B}_i$ is a set of linearly independent vectors in
$M(\underline{v})/M(\underline{v}+\underline{1})$. For $x \in
V_{\{1\} \cup K}$, where $K=\{a_1,\cdots,a_k\} \subset
\{2,\cdots,r\}$ and $k \geq 1$, one sees easily that
$\{{\partial_{k}}|_{V_{\{1\} \cup K \setminus \{a_1\}}}(x), \ldots,
{\partial_{k}}|_{V_{\{1\} \cup K \setminus \{a_k\}}}(x)\}$ is a
linearly independent set. Let $i \leq k$ and $(I_l)_{l \in L}
\subset \{1\} \cup K$, with $|I_l|=i$ for $l \in L$. We denote
$\I=\cup_{l \in L} I_l$.
We have
\begin{eqnarray}
{\partial_{k}}|_{V_K}(x) = \sum_{j=1}^k (-1)^{j+1}
{\partial_{k}}|_{V_{\{1\} \cup K \setminus \{a_j\}}}(x)
\end{eqnarray}
from which it follows that
\begin{eqnarray*}
\textnormal{dim }\left( \bigoplus_{l \in L}
<{\partial_{i}}|_{V_{I_l}}(x)>\right) & = &
\sum_{I_l, 1 \in I_l} 1 
 =  \sum_{I_l, 1 \in I_l} 1 + \left(\sum_{I_l, 1 \notin I_l} 1 -
\sum_{\substack{I \subset \I, |I|=i+1 \\ 1 \in I}} 1 \right) \\
& + & \left(-\sum_{\substack{I \subset \I,|I|=i+1 \\ 1 \notin I}} 1
+
\sum_{\substack{I \subset \I, |I|=i+2 \\ 1 \in I}} 1 \right) \\
& + & \left(\sum_{\substack{I \subset \I,|I|=i+2 \\ 1 \notin I}} 1 -
\sum_{\substack{I \subset \I, |I|=i+3 \\ 1 \in I}} 1 \right)  +  \ldots
\end{eqnarray*}
We can now deduce that
\[\textnormal{dim (Im}
\partial_i)=\textnormal{ dim}C_i-\textnormal{ dim}C_{i+1}+\textnormal{
dim}C_{i+2}-\cdots \qquad \qquad \qquad \qquad \qquad \qquad \blacksquare\] 
\noindent
To prove Corollary \ref{coroltoric} below we need the following lemma.
\begin{lemma}\cite[p. 62]{matsumura} \label{krull}
Let $(A,\mathsf{m})$ be a local Noetherian ring. If $M$ is a finite
$A$-module and $N \subset M$ a submodule, then
\[\bigcap _{n > 0} (N+ \mathsf{m}^n M)=N.\]
\end{lemma}
\noindent The following corollary
reflects well the nature of affine toric varieties in these homological terms.
\begin{corollary} \label{coroltoric}
If $X$ is an affine toric variety and $\nu_1,\ldots,\nu_r$ are monomial valuations, then $h_i^{\underline{v}}=0$ for
$i \geq 1$ and in particular, $P^{\underline{\nu}}_{X}(\underline{t})=\sum_{\underline{v} \in \Z^r} h_0^{\underline{v}}\underline{t}^{\underline{v}}$ is well-defined if one of the following two
equivalent conditions is satisfied:
\begin{enumerate}
 \item \mbox{dim } $\frac{M(\underline{v})}{M(\underline{v}+\underline{e}_1)+\cdots + M(\underline{v}+\underline{e}_r)} < \infty$;
 \item there exists $T \in \Z_{>0}$ such that $M(\underline{v}) \cap \mathsf{m}^T \subset M(\underline{v}+\underline{e}_1)+\cdots + M(\underline{v}+\underline{e}_r)$.
\end{enumerate}
\end{corollary}
\noindent \emph{Proof.} \quad
It follows immediately from Proposition \ref{propbases} that $h_i^{\underline{v}}=0$ for
$i \geq 1$ and obviously $h_0^{\underline{v}}=\mbox{dim } \frac{M(\underline{v})}{M(\underline{v}+\underline{e}_1)+\cdots + M(\underline{v}+\underline{e}_r)}$.
We now show the equivalence between Condition \emph{1} and \emph{2}.
If Condition \emph{2} holds, then  \[\mbox{dim }\frac{M(\underline{v})}{M(\underline{v}+\underline{e}_1)+\cdots + M(\underline{v}+\underline{e}_r)}\leq \mbox{dim } \frac{M(\underline{v})}{M(\underline{v}) \cap \mathsf{m}^T}
\leq \mbox{ dim } \mathcal{O}_{X,o}/\mathsf{m}^T\]
and hence is finite.

Suppose now that Condition \emph{1} holds. Say $\frac{M(\underline{v})}{M(\underline{v}+\underline{e}_1)+\cdots + M(\underline{v}+\underline{e}_r)}$ is $m$-dimensional with basis $\{\overline{f_1},\ldots,\overline{f_m}\}$.
For $T \in \Z_{>0}$, we consider the vector space
\[L_T:=\{(\lambda_1,\ldots,\lambda_m) \in \C^m \mbox{ $|$ } \sum_{i=1}^m \lambda_i f_i \in M(\underline{v}+\underline{e}_1)+\cdots + M(\underline{v}+\underline{e}_r)+(M(\underline{v}) \cap \mathsf{m}^T)\}.\]
We have $\C^m \supseteq L_T \supseteq L_{T+1} \supseteq \ldots$ and hence this chain stabilizes from some $T_0$ on.
Then for $(\lambda_1,\ldots,\lambda_m) \in L_{T_0}$, we have
\begin{eqnarray*}
\sum_{i=1}^m \lambda_i f_i & \in & \bigcap_{T \geq T_0} L_T \\
& \subseteq & \bigcap_{T \geq T_0} \left(M(\underline{v}+\underline{e}_1)+\cdots + M(\underline{v}+\underline{e}_r)+\mathsf{m}^T \right) \\
& \subseteq & \bigcap_{T \geq T_0} \left(M(\underline{v}+\underline{e}_1)+\cdots + M(\underline{v}+\underline{e}_r)\right) \qquad (\mbox{Lemma } \ref{krull}).
\end{eqnarray*}
This implies that $L_{T}=0$ for $T \geq T_0$. It now follows easily for all $T \geq T_0$ that
$M(\underline{v}) \cap \mathsf{m}^T \subset M(\underline{v}+\underline{e}_1)+\cdots + M(\underline{v}+\underline{e}_r)$.
\hfill $\blacksquare$

\begin{remark} \label{opmerkingtorischCd} \emph{If $X= \C^d$ equipped with monomial valuations $\nu_j=(\nu_{j1},\ldots,\nu_{jn})$, $1 \leq j \leq r$, then $P^{\underline{\nu}}_X(\underline{t})$ is well-defined if for all $i \in \{1,\ldots,n\}$, there exists some $j \in \{1,\ldots,r\}$ such that $\nu_{ji} \neq 0$.}
\end{remark}
\noindent The above observations make it natural to propose the following definition for Poincar\'e series, which
we will denote by $\textit{\textbf{P}}^{\underline{\nu}}_X(\underline{t})$.
\begin{definition}
Let $(X,o)$ be a germ of a complex algebraic variety and let $\underline{\nu}=\{\nu_1,\cdots,\nu_r\}$ be a set of
discrete order functions on $\mathcal{O}_{X,o}$. If all $h_i^{\underline{v}}$, defined as above, are finite, then
the Poincar\'e series $\textbf{P}^{\underline{\nu}}_X(\underline{t})$ is defined as
\[\sum_{\underline{v} \in \mathbb{Z}^r} \left(\sum_{i=0}^{r-1}(-1)^i h_i^{\underline{v}}\right)\underline{t}^{\underline{v}}.\]
\end{definition}

${}$
\begin{center}
\textsc{3. Poincar\'e series of embedded filtrations}
\end{center}
${}$ \\
\noindent In \cite{annembedded}, the second author introduced a Poincar\'e series for an embedded subspace $V$ defined by an ideal
$\I$ in $\mathcal{O}_{X,o}$ and order functions $\underline{\nu}=\{\nu_1,\cdots,\nu_r\}$ on $\mathcal{O}_{X,o}$.
Keeping the notation from Section 1, this Poincar\'e series was defined as
\[\mathcal{P}^{\underline{\nu}}_{V}(t_1,\cdots,t_r):=\frac{\prod_{j=1}^{r}(t_j-1)}{(t_1\cdots t_r - 1)}\sum_{\underline{v}
\in \mathbb{Z}^r} \mbox{dim}(M(\underline{v})+\I/M(\underline{v}+\underline{1})+\I)
\underline{t}^{\underline{v}},\]
if the dimension of $(M(\underline{v})+\I/M(\underline{v}+\underline{1})+\I)$ is finite.
Similarly as above, this description can be generalized to a homological description.

Let $V^{\I}_I:=(M(\underline{v}+\underline{e}_I)+\I)/(M(\underline{v}+\underline{1})+\I)$
for $I \subseteq \{1,\ldots,r\}$.
Let \[ C^{\I}_i:=\begin{cases}
\begin{array}{cc}
\{0\} & \mbox{ if } i=-1, i>r, \\
\bigoplus_{\substack{I
\subset \{1,\ldots,r\}, \\ \#I=i}} V^{\I}_I & \mbox{ if } 0 \leq i \leq r.
\end{array}
\end{cases} \]
For $-1 \leq i \leq r-1$, we define a map $\partial^{\I}_{i+1}:
C^{\I}_{i+1} \rightarrow C^{\I}_i$ by defining it on each component $V^{\I}_I$ ($\#
I=i+1$). Suppose that $I=\{a_1,a_2,\ldots,a_{i+1}\}\subseteq \{1,\ldots,r\}$, with $a_1 < a_2< \ldots <a_{i+1}$,
then we set
\begin{eqnarray*}
V^{\I}_I & \longrightarrow & C^{\I}_i=\bigoplus_{\substack{J
\subset \{1,\ldots,r\}, \\{\#J=i}}} V^{\I}_{J} \\
x & \longmapsto & (\underline{y})_J,
\end{eqnarray*}
where $y_J=0$ if $J \nsubseteq I$ and $y_J=(-1)^k x$ if $J = I
\setminus \{a_k\}$. We get a complex
$\mathcal{C}^{\I,\underline{v}}_\bullet=\mathcal{C}^{\I}_\bullet=(C^{\I}_i,\partial^{\I}_i)_{i \in \mathbb{Z}_{\geq -1}}$ of vector spaces and we will denote the arising
homologies here by $h^{\I,\underline{v}}_{i}$. We give the following definition.
\begin{definition}
Let $(X,o)$ be a germ of a complex algebraic variety and let $\underline{\nu}=\{\nu_1,\ldots,\nu_r\}$ be a set of
discrete order functions on $\mathcal{O}_{X,o}$. Let $\I$ be an ideal in $\mathcal{O}_{X,o}$ defining the subspace $V$. If all $h^{\I,\underline{v}}_{i}$ are finite, then
the Poincar\'e series $\mathbfcal{P}^{\underline{\nu}}_V(\underline{t})$ is defined as
\[\sum_{\underline{v} \in \mathbb{Z}^r} \left(\sum_{i=0}^{r-1}(-1)^i h^{\I,\underline{v}}_{i}\right)\underline{t}^{\underline{v}}.\]
\end{definition}
\noindent When the order functions are all valuations and when the ideal $\I=(h)$ is principal, then it was shown in
\cite[Theorem 1]{annembedded} that
\[\mathcal{P}^{\underline{\nu}}_V(\underline{t})=(1-\underline{t}^{\underline{q}})
P^{\underline{\nu}}_X(\underline{t})\]
where $\underline{q}=\underline{\nu}(h)$.
This proof worked whenever at least one valuation was centered in the maximal ideal. We now want to extend this formula for the Poincar\'e series defined in homological terms, i.e.\ including sets of valuations where none of them is centered at the maximal ideal.
\begin{theorem} \label{thmembeddednotcentered}
Let $(X,o)$ be irreducible, $\I=(h)$ a principal ideal in
$\mathcal{O}_{X,o}$ and $V$ the analytic subspace of $(X,o)$ determined by the ideal $\I$. Let $\underline{\nu}=\{\nu_1,\ldots,\nu_r\}$ be a set of discrete valuations of
$\mathbb{C}(X)$ such that $\textbf{P}^{\underline{\nu}}_X(\underline{t})$ and $\mathbfcal{P}^{\underline{\nu}}_V(\underline{t})$ are well-defined.
We write $\underline{q}=\underline{\nu}(h)$. Then
\begin{eqnarray}\label{compunotcentered}
\mathbfcal{P}^{\underline{\nu}}_{V}(\underline{t})=(1-\underline{t}^{\underline{q}})
\textbf{P}^{\underline{\nu}}_X(\underline{t}).
\end{eqnarray}
\end{theorem}
\noindent \emph{Proof.} \quad Let us study the coefficient of $\underline{t}^{\underline{v}}$ in both members.
With the notation of Section 2, the coefficient of $\underline{t}^{\underline{v}}$ in the right hand side member is
$\sum_{i=0}^{r-1}(-1)^i h_i^{\underline{v}}-\sum_{i=0}^{r-1}(-1)^i h_i^{\underline{v}-\underline{q}}$.
With the notation of Section 3, the coefficient of $\underline{t}^{\underline{v}}$ in the left hand side member can be written as $\sum_{i=0}^{r-1}(-1)^i h^{\I,\underline{v}}_{i}$. Now we will make the link between the three complexes that arise. Notice that
\begin{eqnarray*}
0 \rightarrow M(\underline{v}-\underline{q}) \stackrel{\alpha}{\rightarrow} M(\underline{v}) \stackrel{\pi}{\rightarrow} M(\underline{v})/(h)M(\underline{v}-\underline{q}) \rightarrow 0,
\end{eqnarray*}
is an exact sequence, where $\alpha$ is the multiplication map with $h$ and $\pi$ is the projection map.
As $\nu_1,\ldots,\nu_r$ are valuations, one has $(h)M(\underline{v}-\underline{q})=(h) \cap M(\underline{v})$, and so
$M(\underline{v})/(h)M(\underline{v}-\underline{q}) \cong M(\underline{v})+(h)/(h)$.
This makes that we have a short exact sequence of the complexes
\begin{eqnarray*}
0 \rightarrow  \mathcal{C}^{\underline{v}-\underline{q}}_\bullet
\rightarrow \mathcal{C}^{\underline{v}}_\bullet \rightarrow \mathcal{C}^{\I,\underline{v}}_\bullet \rightarrow 0.
\end{eqnarray*}
This sequence induces a long exact sequence of the homology groups (see for example \cite[p. 203]{ha}) from which one then easily deduces that
$\sum_{i=0}^{r-1}(-1)^i h^{\I,\underline{v}}_{i}=
\sum_{i=0}^{r-1}(-1)^i h_i^{\underline{v}}-\sum_{i=0}^{r-1}(-1)^i h_i^{\underline{v}-\underline{q}}$.
\hfill \begin{flushright} $\blacksquare$ \end{flushright}
\begin{remark} \label{remexistencepoincare}
Notice that it follows from the proof of Theorem \ref{thmembeddednotcentered} that the existence of $\textbf{P}^{\underline{\nu}}_X(\underline{t})$ implies the existence of $\mathbfcal{P}^{\underline{\nu}}_V(\underline{t})$ for every analytic subspace $V$ of $(X,o)$ determined by a principal ideal $\I$.
\end{remark}
As this relation between the embedded Poincar\'e series and the ambient Poincar\'e series holds in this broader context, one can extend properties such as \cite[Thm. 6]{annembedded} about functions $h \in \C[x_1,\ldots,x_n]$ that are nondegenerate with respect to their Newton polyhedron $\mathcal{N}$. The condition that $\mathcal{N}$ has at least one compact facet in \cite[Thm. 6]{annembedded} (notice that a compact facet induces a valuation centered at the maximal ideal) can now be weakened to the condition that the variable $x_i$ does not divide the polynomial expression of $h$ if the monomial $x_i$ is already contained in the support of $h$, $1 \leq i \leq d$. The property becomes then:
\begin{theorem}
Suppose that $h \in \C[x_1,\ldots,x_n]$ is nondegenerate with respect to
its Newton polyhedron $\mathcal{N}$ in the origin and that no monomial $x_i$ is dividing the polynomial expression of $h$ for which
$x_i$ is already contained in the support of $h$, $1 \leq i \leq d$.
Let $\underline{\nu}=\{\nu_1,\cdots,\nu_r\}$ be the monomial valuations on $\C^d$ induced by the facets of $\mathcal{N}$.
Then the Poincar\'e series $\mathbfcal{P}^{\underline{\nu}}_{V}(\underline{t})$ contains the same information as the Newton polyhedron of $h$ and in particular determines the zeta function
of monodromy of $h$.
\end{theorem}
\noindent \emph{Proof.} \quad
It follows from Remark \ref{opmerkingtorischCd} that $\textbf{P}^{\underline{\nu}}_{\C^d}(\underline{t})$ is well defined because $\{\nu_1,\cdots,\nu_r\}$ contains the valuations $\nu_1:=(1,0,\ldots,0), \ldots, \nu_d:=(0,\ldots,0,1)$, and hence by Remark \ref{remexistencepoincare} that also $\mathbfcal{P}^{\underline{\nu}}_{V}(\underline{t})$ exists.
With $\underline{q}:=\underline{\nu}(h)$, we know from Theorem \ref{thmembeddednotcentered} that
\begin{eqnarray} \label{eqnrelation}
\mathbfcal{P}^{\underline{\nu}}_{V}(\underline{t})=
(1-\underline{t}^{\underline{q}})\textbf{P}^{\underline{\nu}}_{\C^d}(\underline{t}).
\end{eqnarray}
No factors cancel in Equation (\ref{eqnrelation}). Indeed, suppose that $\underline{\nu}$ contains only the valuations $\nu_1, \ldots, \nu_d$, then the form of the Newton polyhedron implies that $h$ is a constant and so $\underline{q}=\underline{0}$ contradicting the cancelation, or that $h$ is a multiple of some monomial $m$. Then $\underline{\nu}(h)=\underline{\nu}(m)$ but a cancelation would imply that $m=x_i$ for some $i \in \{1,\ldots,d\}$ what contradicts the hypothesis on $h$.
Suppose now that $\underline{\nu}$ contains more valuations, say we have $\nu_{d+1}=(\nu_{d+1,1},\ldots,\nu_{d+1,d})$ with at least two entries different from $0$. This means that $\mathcal{N}$ contains a facet such that the affine space generated by this facet has as equation
\[\nu_{d+1,1}x_1 + \cdots + \nu_{d+1,d}x_d=N,\]
with $N \geq \nu_{d+1,1} +\cdots + \nu_{d+1,d}$. Notice that $N$ is also equal to $\nu_{d+1}(h)=q_{d+1}$. If there would be some cancelation, then $\underline{q}=(\nu_{1,i},\ldots,\nu_{r,i})$ for some $i \in \{1,\ldots,d\}$. Thus $N=q_{d+1}=\nu_{d+1,i}$ which contradicts the fact that $\nu_{d+1}$ contains at least two entries different from $0$.
\hfill $\blacksquare$
${}$ \\ \\
\noindent Antonio Campillo; Universidad de Valladolid, Departamento de
\'Algebra, Geometr\'ia y Topolog\'ia, Valladolid, Espa\~{n}a, email:
campillo@agt.uva.es
\\ \\
Ann Lemahieu; Universit\'e Lille 1, U.F.R. de Math\'ematiques, 59655 Villeneuve d'Ascq C\'edex, email: ann.lemahieu@math.univ-lille1.fr
\\ \\
Acknowledgement: The research was partially supported by the Fund of
Scientific Research - Flanders and MEC PN I+D+I MTM2007-64704.

\footnotesize{
\end{document}